\documentclass[10pt]{article}
\catcode`\@=11 \catcode`\@=12

\usepackage{amsmath}
\usepackage{amssymb}
\usepackage{amsthm}
\usepackage{oldgerm}

\usepackage[utf8]{inputenc}

\usepackage[pdftex,colorlinks,urlcolor=blue,pdfstartview=FitH]{hyperref}

\newcommand{\Rm}{\mathbb{R}}

\newcommand{\mF}{\ensuremath{\mathcal{F}}}

\newcommand{\Nm}{\ensuremath{\mathbb{N}}}

\newcommand{\Tm}{\ensuremath{\mathbb{T}}}

\newtheorem{lem}{Lemma}
\newtheorem{ad}{Addendum}
\newtheorem{hyp}{Hypothesis}
\newtheorem{thm}{Theorem}

\newtheorem{cor}[lem]{Corollary}
\newtheorem{prop}[lem]{Proposition}

\newtheorem{defn}[lem]{Definition}

\def\proof {\noindent{\sc{Proof. }}}
\def\qed {\mbox{}\hfill {\small \fbox{}} \\}
\def\lto{\longrightarrow}
\def\lmto{\longmapsto}

\def\leq{\leqslant}
\def\geq{\geqslant}

\author{P. BERNARD}
\title{Semi-concave singularities and the Hamilton-Jacobi equation}

\begin{document}

\maketitle

\begin{abstract}
We study the Cauchy problem for the Hamilton-Jacobi equation with a semi-concave 
initial condition. We  prove an inequality between the two types
of weak solutions emanating from such an initial condition (the variational and the viscosity solution). We also give conditions for an explicit semi-concave function
to be a viscosity solution. These conditions generalize the entropy inequality
characterizing piecewise smooth solutions of scalar conservation laws in dimension one.

\end{abstract}

\section{Introduction}

We  consider the Cauchy problem for the Hamilton-Jacobi equation
\begin{equation}\tag{HJ}\label{HJ}
\partial_t u(t,x) +H(t,x,\partial_x u(t,x))=0,
\end{equation} 
of an unknown function $u(t,x):\Rm\times \Rm^d\lto \Rm$.
It will also be useful to consider the associated Hamiltonian system
\begin{equation}\label{HS}\tag{HS}
\dot q(t)=-\partial_pH(t,q(t),p(t))\quad,\quad 
\dot p(t)=\partial_qH(t,q(t),p(t)).
\end{equation}
We will most of the time assume :
\begin{hyp}\label{1}
The Hamiltonian $H(t,x,p):\Rm\times \Rm^d\times \Rm^d\lto \Rm$
is $C^2$ and there exists a constant $A$ such that 
$$
|d^2H(t,x,p)|\leq A, \quad |dH(t,x,p)|\leq A(1+|p|), \quad |H(t,x,p)|\leq A (1+|p|)^2
$$
for each $(t,x,p)$.
\end{hyp}
In particular, the Hamiltonian system
is complete. No convexity assumption is made on $H$.

We focus our attention on the case where the initial condition $u_0$ 
 is semi-concave and Lipschitz, given as the infimum of an
equi-Lipschitz  family $F_0$ of $C^2$ functions
with uniformly equi-bounded second derivatives, which means that 
there exists a constant $B$ such that
$|d^2f_0(x)|\leq B$ for each $x\in \Rm^d, f_0\in F_0$.
See Section \ref{secnonsmooth} for more details on semi-concave functions.
The general theory of Hamilton Jacobi equations 
allows to define two solutions for the Cauchy problem 
with the Lipschitz initial condition $u_0$ at time $0$:

The  variational solution $g(t,x)=G_0^tu_0(x)$, and
the viscosity solution $v(t,x)=V_0^tu_0(x)$, see Section \ref{secCauchy}
for more details.
One of our goals in the present work is to compare these two solutions
and the following natural third candidate.

The theory of characteristics implies that there exists a constant $T(B)>0$,
which depends on $B$ (and on $A$), such that, to each function $f_0\in F_0$
is associated a $C^2$ solution $f:]-T(B),T(B)[\times \Rm^d\lto \Rm$ of (\ref{HJ})
satisfying $f(0,.)=f_0$, more details in Section \ref{secCauchy}.
We denote by $F$ this family of $C^2$ solutions. 
Their infimum is a natural candidate to be a solution of our Cauchy problem
on $[0,T(B)[\times \Rm^d$, although it depends on the family $F_0$, and not just 
on the function $u_0$.

\begin{thm}\label{thm1}
The solutions $v$ and $g$ are semi-concave on $[0,T[\times \Rm^d$
for each $T\in ]0,T(B)[$, and 
$$
v\leq g\leq \inf _{f\in F} f
$$
on $[0, T(B)[\times \Rm^d$.
\end{thm}

This theorem is proved in section \ref{proofthm1}, where a sufficient condition
for the equality $g=\inf f$ is also given.
We discuss this condition here under the additional assumption that $F$ is 
closed for the $C^1_{loc}$ topology
(this is a very minor restriction since one can always replace $F$ by its closure).
In this case, we denote
by $\partial_xF(t,x)$ the set $\{\partial_x \underline f(t,x), \underline f\in F, 
\underline f(t,x)=\min f (t,x)\}$.

\begin{ad}
If $\partial_x F(0,x)$ is convex for each $x$, then the equality $g=\inf f$
holds in Theorem \ref{thm1}.

If  $\partial_x F(t,x)$ is convex for each $(t,x)$, then the equality $v=g=\inf f$
holds.
\end{ad}

The sufficient condition for the equality $v=\inf f$ mentioned above is 
actually too demanding,
for example it is usually not satisfied in the context of the Hopf formula
for concave solutions, see below. We now propose more reasonable sufficient conditions
inspired by the famous entropy inequalities which characterize 
piecewise smooth entropy solutions of conservation laws, see \cite{Lax}.

We  consider a semi-concave function $u(t,x):]0,T[\times \Rm^d \lto \Rm$ 
which solves (\ref{HJ}) at each point of differentiability.
It can be for example  the variational solution $g$ or the function $\inf f$
in Theorem \ref{thm1}.
We  call such a function a semi-concave solution.
We denote by $D^eu(t,x)$ the set of extremal points of the 
super-differential $Du(t,x)$ of $u$, and by $D^e_xu(t,x)$ the projection 
of $D^eu(t,x)$ on the spacial directions.
We also denote by $D_xu(t,x)$ the projection on the spacial
directions of $Du(t,x)$. Note that $D^e_xu(t,x)$
is bigger than the set $D^eu_t(x)$ of extremal super-differentials 
of the function $u_t=u(t,.)$.
We then denote by $\check H^u_{t,x}$ the greatest convex function of $\Rm^d$
which is smaller than or equal to $H_{t,x}$ on $D^e_xu(t,x)$.
Similarly, we denote by 
 $\hat H^u_{t,x}$ the smallest concave function of $\Rm^d$
which is greater than or equal to $H_{t,x}$ on $D^e_xu(t,x)$.
These functions take the value $+\infty$ (resp $-\infty$) outside of $D_xu(t,x)$.
The following result  holds for all continuous Hamiltonians (not necessarily satisfying
Hypothesis \ref{1}) :

\begin{thm}\label{entropythm}
Let $u(t,x):]0,T[\times \Rm^d \lto \Rm$  be a semi-concave solution of  (\ref{HJ}). 
If 
\begin{equation}\label{convex}
H_{t,x}\leq \check H ^u_{t,x} \quad \text{on}\quad D_xu(t,x)
\end{equation}
for each $(t,x)\in ]0,T[\times \Rm^d$,
then  $u$ is a viscosity solution of (\ref{HJ}).
Conversely, if $u$ is a viscosity solution of (\ref{HJ})
then the inequality
\begin{equation}\label{concave}
H_{t,x}\leq \hat H ^u_{t,x} \quad \text{on}\quad D_xu(t,x)
\end{equation}
holds for each $(t,x)$.
\end{thm}

Conditions of the same kind were introduced in \cite{MC} in a different setting.
In the case where $u$ is the minimum  of two functions $f^-$ and $f^+$,
conditions (\ref{convex}) and  (\ref{concave}) are equivalent,
and they can be expressed as
\begin{equation}\label{linear}
H(t,x, sp^-+(1-s)p^+)\leq s H(t,x,p^-)+ (1-s) H(t,x,p^+)\quad \forall s\in [0,1],
\end{equation}
where $p^{\pm}=\partial_xf^{\pm}(t,x)$.
Theorem \ref{entropythm} then reads 
\begin{cor}\label{cor}
Let $f^-,f^+:]0,T[\times \Rm^d\lto \Rm$ be two $C^2$ solutions
of the Hamilton-Jacobi equation with bounded second derivative and bounded derivative. 
The  function $u:= \min (f^-,f^+) $
is a viscosity solution of (\ref{HJ}) if and only if the entropy condition 
(\ref{linear}) is satisfied at each point $(t,x)\in ]0,T[\times \Rm^d$.
\end{cor}

In the case $d=1$ this corollary is the counterpart of a standard result 
concerning  the conservation law
\begin{equation}\label{CL}\tag{CL}
\partial_t p(t,x)+\partial_x\big(H(t,x,p(t,x))\big) =0,
\end{equation}
which, formally, is the equation solved by the differential $p=\partial_x u$
of solutions of (\ref{HJ}).
For such equations, there is a theory of entropy solutions, which is the counterpart
of the theory of viscosity solutions.
Let us consider  a solution $p(t,x)$  of (\ref{CL})
which is composed of two smooth
branches of solutions  $p^-(t,x)$ and $p^+(t,x)$ on both sides of
 a discontinuity $\chi(t)$,
with $p(t,x)=p^-(t,x)$ for $x\leq \chi(t)$ and $p=p^+$ for $x\geq \chi(t)$.
Assuming in addition that $p^-\geq p^+$ (a condition satisfied by the derivative
of a semi-concave singularity), it is known that the function $p(t,x)$
is an entropy solution of (\ref{CL}) if and only if the entropy condition 
(\ref{linear}) holds at  $(t,\chi(t))$ for each $t$, see \cite{Lax}.
Corollary \ref{cor} is the transposition of this celebrated result
for viscosity solutions of the Hamilton-Jacobi equation.

All the results presented in this note  have obvious counterparts 
 for semi-convex initial 
conditions. The reason why we preferred to work with semi-concave solutions 
is that they play a special role in the case of convex Hamiltonians
(meaning that $\partial^2_{p}H>0$).
As is well-known, viscosity solutions are variational in this case,
forming a single notion of weak solution $(v=g)$. 
This solution is given by an explicit expression, the so-called 
Lax-Oleinik semi-group, and it has  the property of being locally  semi-concave
whatever the initial condition, see \cite{CaSi} for example.
In this convex case, the generalized entropy condition (\ref{convex}) (hence (\ref{concave}))
always holds, hence each semi-concave solution is a viscosity solution
(this property is well-known, see \cite{CaSi} for example).
As a consequence, equality always holds in Theorem \ref{thm1}.
This could also be proved easily with the Lax-Oleinik formula. 

Still in the case of a convex Hamiltonian, semi-convex solutions
are also of interest. 
The first conclusion of Theorem \ref{1} (the solution emanating from
a semi-convex initial data is semi-convex) was first proved in \cite{ENS}.
This solution is then locally $C^{1,1}$ (because a function 
which is both semi-convex and semi-concave is $C^{1,1}$).
This is reflected in Theorem \ref{entropythm} as follows:
The necessary condition to be a semi-convex viscosity solution
 reads $H_{t,x}\geq \check H^u_{t,x}$. 
In view of the assumption of strict convexity of $H$, this can hold
only if the sub-differential $D_xu(t,x)$ is reduced to a point, which implies that the function $u$
is $C^1$.

All these considerations  indicate that Theorem \ref{entropythm}
is not useful in the convex case. To illustrate its possible usefulness,
let us use it to recover the formula of Hopf for concave solutions of integrable Hamiltonians.
We suppose in this discussion that $H=H(p)$ does not explicitly depend on $t$ and $x$,
and consider a concave initial condition $u_0$, that we also assume Lipschitz
for simplicity (this assumption is removed in Section \ref{secHF}).
We write $u_0$ as an infimum  of affine
functions  with the formula
$
u_0(x)=\inf_p( p\cdot x -u_0^*(p)),
$
 where $u_0^*(p)$ is the (concave version of the) Legendre transform
$u_0^*(p):=\inf_x (p\cdot x-u_0(x))$. 
This dual  $u^*_0$  takes the value $-\infty$, outside of a bounded domain 
denoted by $P$.
The 
function $f(t,x)=px-u_0^*(p)-tH(p)$ is a
$C^2$ solution emanating from the affine initial 
condition $f_0(x)=px-u_0^*(p)$.

\begin{cor}\label{corhopf}
If $H(p)$ is a continuous Hamiltonian, and $u_0$ is a Lipschitz and concave 
initial condition, then the function
$$u(t,x):=\min_{p\in P} \big( px-u_0^*(p)-tH(p)\big)
$$ 
is a viscosity solution on $[0,\infty[\times \Rm^d$.
If, in addition, $H$ satisfies Hypothesis \ref{1}, then
$$v(t,x)=g(t,x)=\min_{p\in P} \big( px-u_0^*(p)-tH(p)\big).$$
\end{cor}

That the right hand side in this expression is actually a viscosity 
solution is well-known, even in broader contexts, see for example \cite{LiR,imbert}.

The paper is organized as follows:
In Section \ref{secCauchy}, we quickly recall some basic facts on the various
notions of solutions of the Cauchy problem.
In Section \ref{secnonsmooth}, we settle some notations and elementary properties
on semi-concave functions seen as infima of $C^2$ functions.
We then prove Theorem \ref{thm1} and its addendum in Section \ref{proofthm1}, using a Proposition 
proved in \cite{ENS} in the case of a convex Hamiltonian. 
We prove Theorem \ref{entropythm} in Section \ref{secentropy}, in the general setting 
of a contiuous Hamiltonian. We return to the Hopf formula and prove Corollary \ref{corhopf}
in Section \ref{secHF}.

\section{The Cauchy Problem}\label{secCauchy}
We give here a very brief survey
on  the  Cauchy problem associated to  the
Hamilton-Jacobi equation (\ref{HJ}).

\subsection{Classical solutions}\label{classical}

The theory of characteristics links classical solutions
of (\ref{HJ}) with the Hamiltonian system (\ref{HS}). 
We give here a brief account on the results, see for example \cite{CaSi,LO} for more details.

For each $C^2$ initial condition $f_s$ with bounded second derivative, 
there exist a time  $T>0$ and a $C^2$ solution  $f:]s-T,s+T[\times \Rm^d\lto \Rm$
of (\ref{HJ}) satisfying $f(s,.)=f_s$.
It is necessary to be more quantitative.

There exists a non-decreasing semi-group $Q^t(r)$ on $[0,\infty]$,
such that the time  $T(r):= \sup\{t\geq 0, Q^t(r)<\infty\}$
is positive for each $r\in [0,\infty[$ and such that:

For each $C^2$ initial data $f_s$ satisfying 
$\|d^2f_s\|_{\infty}\leq r$, there exists a unique $C^2$
solution $f:]s-T(r),s+T(r)[\times \Rm^d\lto \Rm$. This solution satisfies
$$
\|d^2f_t\|_{\infty}\leq Q^{|t-s|}(r).
$$
If moreover $f_s$ is Lipschitz (that is, if $df_s$ is bounded),
we have the estimates 
\begin{itemize}
\item
$Lip(f_t)\leq (Lip(f_s)+1)e^{A|t-s|}-1$
\item
 $|\partial_tf(t,x)|\leq A\big((Lip(f_s)+1)e^{A|t-s|}\big)^2$
\item
$|d^2f(t,x)|\leq D \quad \forall (t,x)\in ]s-T,s+T[\times \Rm^d,$
\end{itemize}
for each $T<T(r)$, where $D$ is some constant depending only on $A,Lip(f_s),Q^T(B)$.
 
This solution is related to the Hamiltonian system as follows:
For each $(t,x)\in ]-T,T[\times \Rm^d$, the Hamiltonian trajectory 
$(q(s),p(s))$ which satisfies 
$(q(t),p(t))=(x,\partial_xf(t,x))$ satisfies
$$
p(s)=\partial_xf(s,q(s))
$$
for each $s\in ]-T,T[$, and
$$
f(t,x)=f(s,q(s))+\int_s^t p(s)\cdot \dot q(s) -H(s,q(s),p(s)) ds.
$$
It is usually not possible to extend $C^2$ solutions to the whole real line,
which led to the introduction of some notions of weak solution.

\subsection{Variational solutions}

See for example \cite{BC,Vit} for more on variational solutions.
This notion of solutions is directly inspired by the method of characteristics.

A function $g(t,x):[s,\infty)\times \Rm^d\lto \Rm$ is called a variational solution of the Cauchy problem
with Lipschitz initial data $u(x)$ at time $s$ if it is locally Lipschitz,
satisfies the initial condition $g_s=u$,
solves the equation almost everywhere, and
if, for each $(t,x)\in ]s,\infty)\times \Rm^d$, there exists
a trajectory $(q(s),p(s))$ of the Hamiltonian system such that 
$q(t)=x$, $p(0)\in Du(q(0))$, and
$$
u(t,x)=u_0(q(0))+\int_0^t p(s)\cdot \dot q(s)-H(s,q(s),p(s)) ds.
$$
Here $Du$ denotes the Clarke differential of $u$, see Section \ref{secnonsmooth}.
In other words, $g(t,x)$ is a critical value of the functional
$$
(q(s),p(s))\lmto u(q(0))+\int_0^t p(s)\cdot \dot q(s)-H(s,q(s),p(s)) ds
$$
on the space of $C^1$ curves 
$(q(s),p(s)):[0,t]\lto \Rm^d\times \Rm^d$ 
which satisfy 
$q(t)=x.
$

There exists a family of operators $G_s^t$, $s\leq t$
which map $C^{0,1}(\Rm^d)$ (the space of Lipschitz functions)
into itself,
such that the function  $(t,x) \lmto G_s^tu (x)$ is a variational solution
with initial data $u$ at time $s$,
 and such that 
\begin{enumerate}
\item $Lip(G_s^tu)\leq (Lip(u)+1)e^{A(t-s)}-1$,\quad
 $\|G_s^tu-u\|_{\infty}\leq A(t-s)\big((Lip(u)+1)e^{A(t-s)}\big)^2$
\item $u\leq v\Rightarrow G_s^tu\leq G_s^tv$
\item If $f(t,x):]T^-,T^+[\times \Rm^d \lto \Rm$ is a $C^2$ solution, 
then $G_s^t f_s=f_t$ for each $s\leq t$ in $]T^-,T^+[$.
\end{enumerate}

A family of operators $G_s^t$ satisfying the properties 
above is called a variational resolution of (\ref{HJ}).
There is no uniqueness for variational solutions, 
and not even uniqueness for variational resolutions.
It would be tempting to ask in addition that the resolution $G_s^t$
satisfy the Markov property 
$G_s^t\circ G_{\tau}^s=G_{\tau}^t$.
However, adding such a condition to the properties $(1-3)$ above would
lead to the Viscosity resolution, see below, which does not 
produce variational solutions in general.

We consider that a variational resolution $G$ is fixed once and for all
in the present paper. When we speak of \textit{the} variational solution
emanating from an initial condition $u_0$, we mean the function 
$g(t,x)=G_0^tu_0(x)$.

\subsection{Viscosity solutions}\label{viscosity}
See for example \cite{CaSi} for more on viscosity solutions.
There exists a unique family of operators $V_s^t$ acting on 
$C^{0,1}(\Rm^d)$ and such that 
\begin{enumerate}
\item $Lip(V_s^tu)\leq (Lip(u)+1)e^{A(t-s)}-1$,\quad
 $\|V_s^tu-u\|_{\infty}\leq A(t-s)\big((Lip(u)+1)e^{A(t-s)}\big)^2$
\item $u\leq v\Rightarrow V_s^tu\leq V_s^tv$
\item If $f(t,x):]T^-,T^+[\times \Rm^d \lto \Rm$ is a $C^2$ solution, 
then $V_s^t f_s=f_t$ for each $s\leq t$ in $]T^-,T^+[$.
\item $V_s^t\circ V_{\tau}^s=V_{\tau}^t$ for each $\tau\leq t\leq s$.
\end{enumerate}

For each Lipschitz initial condition $u$, and each initial time 
$s$, the functions 
$$
]s,\infty[\times \Rm^d \ni (t,x)\lmto V_s^t u(x)
$$
is a viscosity solution of (\ref{HJ}) in the classical sense.
It means that each smooth function $\phi(t,x)$ which
has a contact from above (resp. from below)  with $u$ at some point 
$(t_0,x_0)\in ]s,\infty[\times \Rm^d$ must satisfy 
$$\partial_t \phi(t_0,x_0) +H(t_0,x_0,\partial_x\phi(t_0,x_0))\leq 0,
\quad (\text{resp.} \geq 0).
$$
Conversely, the function $(t,x)\lmto V_s^t u(x)$ is the unique viscosity 
solution $v$ such that $v(s,.)=u$ and such that $v(t,x)-u(x)$ is bounded
and Lipschitz on $]s,T[\times \Rm^d$ for each $T>s$.

The theory of viscosity solutions thus provides a "good" resolution
of the Cauchy problem (existence and uniqueness).
As we mentioned above, apart in some special cases (for example when $H$
is convex in $p$), the viscosity solution is not 
in general a variational solution.

It is an interesting problem in general to describe and compare these
two notions of solutions. 
Let us mention in this direction  a recent statement recently proved by Qiaoling Wei
in \cite{Wei} (see also \cite{Roos}). This result had been conjectured by Chaperon and Viterbo.
We denote, for each $k\in \Nm$,
by $^kG_s^t$ the operator
$$
^kG_s^t:= G^{t}_{t_{[kt]}}\circ G_{t_{[kt]-1}}^{t_{[kt]}}
\circ \cdots \circ G_{t_1}^{t_2}\circ G_{s}^{t_1}
$$
where $t_i=s+i/k$ and $[kt]$ is the integer part of $kt$.
We have 
$$
^kG_s^tu \lto V_s^t u 
$$
locally uniformly (in $t$ and $x$) for each Lipschitz function $u$.

\section{Nonsmooth Calculus and semi-concave functions as minima of $C^2$ functions}\label{secnonsmooth}
We recall some standard definitions and properties of non-smooth 
calculus, see \cite{CaSi}, Chapter 3 and \cite{Clarke}. 
We will consider only locally Lipschitz functions $u:\Rm^d\lto \Rm$.

A super-differential of $u$ at $x$ is a vector $p\in \Rm^d$
such that there exists a $C^1$ function $f$ satisfying $df(x)=p$
and having a contact from above with $u$ at $x$, which means that 
$$
f(x)=u(x)\quad, \quad  f\geq u.
$$

A proximal 
super-differential of $u$ at $x$ is a vector $p\in \Rm^d$
such that there exists a $C^2$ (or, equivalently, smooth) function $f$ satisfying $df(x)=p$
and having a contact from above with $u$ at $x$.
A proximal super-gradient is obviously a super-gradient, but the converse is
not true in general.

For concave functions however, super-differentials in the present
sense coincide with super-differentials in the sense of convex analysis
(the slopes of affine functions which have a contact from above with $u$ at $x$),
hence with proximal super-differentials.

The vector $p\in \Rm^d$ is called a reachable gradient of $u$ at 
$x$ if there exists a sequence $x_n\lto x$ of points of differentiability
of $u$ such that $du(x_n)\lto p$. Recall that the points of differentiability 
of $u$ have full measure (hence it is dense).
The set $D^*u(x)$ of reachable gradients is compact and not empty.

The Clarke differential of $u$, denoted by $Du(x)$ is the convex 
hull, in $\Rm^d$, of $D^*u(x)$, see \cite{Clarke}, Section 2.5.
We denote by $D^eu(x)$ the set of extremal points of $Du(x)$,
note that $D^eu(x)\subset D^*u(x)$.

For a function $u(t,x)$ of two variables, we use the notation
$$
D_xu(t,x):=\{p, \exists \eta, (\eta,p)\in Du(t,x)\}
$$
and similarly for $D^e_xu(t,x)$.

The function $u$ is called semi-concave if $u(x)-B\|x\|^2/2$
is concave for some $B$ (we then say that $u$ is $B$-semi-concave).
For semi-concave functions, the set of super-differentials, 
of proximal super-differentials, and of Clarke differentials coincide.

If $p$ is a super-differential at $x_0$ of the $B$-semi-concave function $u$,
then the function $u(x_0)+p\cdot(x-x_0)+B\|x-x_0\|^2/2$ has a contact from above
with $u$ at $x_0$.

Let us now consider a set  $F\subset C^2(\Rm^d,\Rm)$ with uniformly equi-bounded second derivative
and assume  that the function
$u(x):= \inf _{f\in F} f(x)
$
takes finite values.

The function $u$ is then semi-concave hence  locally Lipschitz.

\begin{defn}\label{dF}
We denote by $dF(x)$ the set of limits of sequences of the form
$df_n(x), f_n\in F, f_n(x)\lto u(x)$. 
\end{defn}

Note that $dF(x)$ is compact.
In the case where $F$ is closed in  $C^1_{loc}$, this is just
$\{df(x), f\in F, f(x)=u(x)\}$.
The set $dF(x)$ depends on $F$, and not only on the function $u$, but it 
is related to the super-differential $Du(x)$:

\begin{lem}\label{diff}
We have $dF(x)\subset Du(x)$ and $dF(x)$ is not empty. 
If $u$ is differentiable at $x$, then $dF(x)=\{du(x)\}$.
\end{lem}

\proof
Let $p=\lim df_n(x)$, with $f_n(x)\lto u(x)$, be a point of 
$dF(x)$.
Since 
$$
u(x+y)\leq f_n(x+y)\leq f_n(x)+df_n(x)\cdot (y-x)+C\|y-x\|^2,
$$
we conclude at the limit that 
$$
u(x+y) \leq u(x)+p\cdot (y-x)+C\|y-x\|^2
$$
hence $p$ is a proximal super-differential of $u$ at $x$, 
$p\in Du(x)$.
We have proved that $dF(x)\subset Du(x)$.

To prove that $dF(x)$ is not empty, we consider
a sequence $f_n\in F$ such that $f_n(x)\lto f(x)$.
We have, for each $y\in \Rm^d$,
$$
u(x+y)\leq f_n(x+y)\leq f_n(x)+df_n(x)\cdot y+C\|y\|^2. 
$$
Applying this inequality with $y_n=-df_n(x)/\|df_n(x)\|$
yields
$$
\|df_n(x)\|\leq f_n(x)-u(x+y_n) +C
$$
which implies that the sequence $df_n(x)$ is bounded.
As a consequence, it has converging subsequences, hence $dF(x)$ is not empty.

If $u$ is differentiable at $x$, then $Du(x)=\{du(x)\}$, hence 
also $dF(x)=\{du(x)\}$.
\qed

\begin{lem}\label{surdiff}
We have $D^*u(x)\subset dF(x)\subset Du(x)$.
As a consequence, the super-differential $Du(x)$ is the convex hull of $dF(x)$.
\end{lem}

\proof
let $p=\lim du(x_n)$ be a reachable gradient, where $x_n\lto x$ is a sequence
of points of differentiability of $u$.
For each $n$, we have $dF(x_n)=\{du(x_n)\}$ hence there exists a function
$f_n\in F$ such that $f_n(x_n)\leq u(x_n)+1/n$ and 
$\|df_n(x_n)-du(x_n)\|\leq 1/n$.
Since $\|df_n(x_n)-df_n(x)\|\leq C\|x_n-x\|$, we conclude that
$df_n(x)\lto p$.
On the other hand we can estimate
\begin{align*}
f_n(x)-u(x)&\leq 
f_n(x_n)+df_n(x_n)\cdot (x-x_n)+C\|x-x_n\|^2-u(x)\\
& \leq df_n(x_n)\cdot (x-x_n)+C\|x-x_n\|^2+f_n(x_n)-u(x_n)+u(x_n)-u(x)
\end{align*}
and the right hand side is converging to $0$.
We conclude that $p\in dF(x)$. 
\qed

\begin{lem}\label{F0}
If $u$ is a B-semi-concave and $L$-Lipschitz function, there exists a 
set  $F\subset C^2(\Rm^d,\Rm)$ such that $|d^2f(x)|\leq B$ 
and $|df(x)|\leq 6L$
for each $f\in F$ and $x\in \Rm^d$, such that 
$u(x)= \min _{f\in F} f(x)
$
and such that $Du(x)=dF(x)$ at each point $x$.
\end{lem}

\proof
Let $\psi(r):[0,\infty)\lto [0,\infty)$ be a non-negative, non-increasing 
function which is equal to $B$ on $[0,4L/B]$ and to $0$ on $[5L/B,\infty)$.
Let $\Psi(r)$ be the primitive of $\psi$ such that $\Psi(0)=0$.
Note that $\Psi(r)\in [0,5L]$ for each $r\in [0,\infty)$.
Let then $\varphi(r)$ be the primitive of $\Psi$ such that $\varphi(0)=0$.
The function $\varphi$ is $5L$-Lipschitz, convex, and it satisfies $0\leq \varphi''\leq B$.
Note also that 
$$
\varphi(r)\geq \min(Br^2/2,2Lr).
$$
Let us consider the family 
$F$ formed by the functions 
$$
u(x_0)+p\cdot (x-x_0)+\varphi(|x-x_0|)
$$
for $x_0\in \Rm^d$ and $p\in Du(x_0)$.
Since we have $|p|\leq L$, these functions are $6L$-Lipschitz.
They also satisfy $|d^2f|\leq B$, see Appendix \ref{hessian}.
It is clear that $Du(x_0)\subset dF(x_0)$ for each $x_0$,
we thus have equality provided that $u=\min_{f\in F} F$,
which is the last thing we have to prove.
It is enough to observe that $f\geq u$ for each $f\in F$.
For all $x_0\in \Rm^d$, $p\in Du(x_0)$ and $x\in \Rm^d$, we have
$$
u(x)\leq u(x_0)+p\cdot (x-x_0)+B|x-x_0|^2
$$
and, since $|p|\leq L$
$$
u(x)\leq u(x_0)+L|x-x_0|\leq u(x_0)+p\cdot (x-x_0)+2L|x-x_0|,
$$
 hence 
$$u(x)\leq u(x_0)+p\cdot (x-x_0)+\varphi(|x-x_0|).
$$
\qed

\section{Proof of Theorem \ref{thm1}.}\label{proofthm1}

We consider a Lipschitz and $B$-semi-concave initial condition $u_0$.

We write this initial condition as
$u_0=\inf_{f_0\in F_0} f_0$ 
for an equi-Lipschitz  family $F_0$ such that 
$|d^2f_0(x)|\leq B$ for all $f_0\in F_0$ and $x\in \Rm^d$.
All the families $F_0$ considered in this section are assumed to satisfy these
conditions.
The initial condition $u_0$ (but not the family $F_0$) and the constant $B$ 
are fixed once and for all.

We  define the family 
 $F\subset C^2([0,T(B)[\times \Rm^d)$
of  solutions of (\ref{HJ}) emanating from elements of $F_0$.
We recall that $F$ is equi-Lipschitz and has uniformly equi-bounded second derivative
on $]0,T[\times \Rm^d$ for each $T\in ]0,T(B)[$.

The inequalities 
$$
G_0^tu_0\leq \inf_{f\in F}f_t
\quad,\quad
V_0^tu_0\leq \inf_{f\in F}f_t
$$
follow from the monotony of the operators $V_0^t$ and $G_0^t$
since $G_0^tf_0=V_0^t f_0=f_t$ for each $f_0\in F_0$.

The key step in the proof of Theorem \ref{thm1} is the following result, which extends 
to the non-convex setting the main 
Proposition of \cite{ENS}:

\begin{prop}\label{hopf}
If $d F_0(x)=Du_0(x)$ for each $x\in \Rm^d$, then 
the inequality $g\geq \inf f$ holds on $[0,T[\times \Rm^d$ for each variational solution $g$. As a consequence, the equality
$$G_0^t u_0=\inf_{f\in F}f_t
$$
holds for each $t\in [0,T[$.
\end{prop}

\proof
Let $g(t,x)$ be a variational solution emanating from $u_0$.
For each $(t_0,x_0)\in ]0,T[\times \Tm^d$, there exists a Hamiltonian orbit  
$(q(s), p(s))$ such that $q(t_0)=x_0$, $p(0)\in D_x u_0(q(0))$,
and 
$$
g(t_0,x_0)=u_0(q(0))+\int_0^{t_0} p(s)\cdot \dot q(s)-H(s,q(s),p(s)) ds. 
$$
Since $p(0)\in D u_0(q(0))$, our hypothesis implies that there exists a
sequence
$f_n\in F$ such that $f_n(0,q(0))\lto u_0(q(0))$ and 
$\partial_xf_n(0,q(0))\lto p(0)$.
Let $(q_n(s),p_n(s)$ be the Hamiltonian trajectory such 
that $q_n(0)=q(0)$ and $p_n(0)=\partial_x f_n (0,q(0)).$
The method of characteristics yields
$$
f_n(t_0,q_n(t_0))=f_n(0,q_n(0))+\int_0^{t_0} p_n(s)\cdot \dot q_n(s)-H(s,q_n(s),p_n(s)) ds.
$$
At the limit, we obtain,
$$
\lim f_n(t_0,x)=u_0(q(0))+\int _0^{t_0} p(s) \cdot \dot q(s)-H(s,q(s),p(s)) ds =g(t_0,x_0).
$$
We conclude that $g\geq \inf f$.
\qed

When the (semi-concave and Lipschitz) initial condition  $u_0$
is given, it is possible to chose the family $F_0$ in such a way that 
the hypothesis  $d F_0(x)=Du_0(x)$ holds for this family, see Lemma \ref{F0}.
We denote by $\mF_0$ a given family with this property, and by $\mF$
the family of associated solutions of $(\ref{HJ})$.
For this family $\mF$, we have the equality 
$g=\inf_{f\in \mF}f$ on $ [0,T(B)[\times \Rm^d$, this is the conclusion of the Proposition we just proved.
Since on the other hand the inequality $v\leq \inf f$
holds for each family $F$, we conclude that 
$$
v \leq g =\inf_{f\in \mF}f\leq \inf_{f\in F}f
$$
on $]0,T(B)[\times \Rm^d$, for each family $F$. 
The equality 
$g=\inf_{f\in \mF}f$ also implies that $g_t$ is semi-concave
with constant $Q^t(B)$(see section \ref{classical} for the definition of $Q^t(B)$)
 and that $g$ is semi-concave on $]0,T[\times \Rm^d$
for each $T\in ]0,T(B)[$.

To prove that the viscosity solution  $v$ is semi-concave, we use the result of Qiaoling Wei, see Section \ref{viscosity},
which states that $^kG_0^tu_0\lto V_0^tu_0$.
Let us denote by $L_0$ the Lipschitz constant of $u_0$, 
and as usual by $B$ its constant of semi-concavity.

For $t\in [0,1/k]\cap ]0,T(B)[$, the function $^kG_0^tu_0=G_0^tu_0$ is $Q^t(B)$-semi-concave
and $L_t$-Lipschitz,
with 
$$
L_t=(L_0+1)e^At-1.
$$
Then, for $t\in [1/k,2/k]\cap ]0,T(B)[$, the function  
$^kG_0^tu_0=G_{1/k}^t\big(G_0^{1/k}u_0\big)$ 
is 
semi-concave with constant $Q^{(t-1/k)}(Q^{1/k}(B))=Q^t(B)$,
and Lipschitz with constant 
$
(L_{1/k}+1)e^{A(t-1/k)}-1=L_t.
$
We prove similarly by recurrence that the function
$^kG_0^tu_0$ is $Q^t(B)$-semi-concave and $L_t$-Lipschitz for each $t\in [0,T(B)[$,
independantly of $k$.
We will use the following Corollary of Proposition \ref{hopf}:

\begin{cor}\label{decroit}
Let $s_0\leq s_1\leq s_2$ be times in $]0,T(B)[$.
Let $u_{s_0}$ be semi-concave with constant $Q^{s_0}(B)$ and
Lipschitz.
Then 
$$
G_{s_1}^{s_2}\circ G_{s_0}^{s_1}u_{s_0}\leq G_{s_0}^{s_2}u_{s_0}.
$$
\end{cor}

\proof
We choose as above, using Lemma \ref{F0} and Proposition \ref{hopf},
a family $\mF_{s_0}$ of $C^2$ functions
such that 
$G_{s_0}^{t}u_{s_0}=\inf _{f_{t}\in \mF_{t}} f_{t}
$
for each $t\in [s_0,T(B)[$.
We have
$$
G_{s_1}^{s_2}\circ G_{s_0}^{s_1}u_{s_0}=
G_{s_1}^{s_2}\big(\inf_{f_{s_1}\in \mF_{s_1}} f_{s_1}\big)
\leq \inf _{f_{s_2}\in \mF_{s_2}} f_{s_2}=G_{s_0}^{s_2}u_{s_0}.
$$
\qed
Let us now fix $T\in ]0,T(B)[$.
Let  $D$ be a constant such that the variational solution emanating from 
a $Q^s(B)$-semi-concave and $L_T$-Lipschitz initial condition $u_s$ 
is $D$-semi-concave on $]s,T[\times \Rm^d$ for each $s\in ]0,T[$.
See Section \ref{classical} for the existence of $D$.
We claim that the function $g^k(t,x):=\, ^kG_0^tu_0(x)$ is $D$-semi-concave for each $k$.

This is true on each interval of the form $]i/k,(i+1)/k[$ where 
$^kg(t,x)=G_{i/k}^t(^kg_{i/k})(x)$.
In addition, the function $^kg$ is equal to
the 
function $G_{i/k}^t(^kg_{i/k})(x)$, on $[i/k,(i+1)/k]$, and, by Corollary
\ref{decroit}, it is not larger than this function on 
$[(i+1)/k,T(B)[$.
Since the function $^kg$ can be touched from above by 
$D$-semi-concave functions at each point, it is $D$-semi-concave.

We conclude, using the result of Qiaoling Wei, that $V_0^tu_0$ is $Q^t(B)$-semi-concave
for each $t\in [0,T(B)[$, and that $v$ is $D$-semiconcave on $[0,T[\times \Rm^d$.
\qed
We now prove the Addendum.
Since $Du_0(x)$ is the convex hull of $dF_0(x)$, the hypothesis
of Proposition \ref{hopf} holds if and only if $dF_0(x)$ is convex.
The first statement of the Addendum is thus a restatement of Proposition \ref{hopf}.

If the hypothesis of the second statement is satisfied, then
$^kG_0^tu_0=u_t$ for each $t$. Using the result of Qiaoling Wei, 
we conclude at the limit that $V_0^tu_0=u_t$ for each $t$.
\qed
\section{Generalized entropy inequalities}\label{secentropy}

We prove Theorem \ref{entropythm}.
The Hamiltonian $H$ is only assumed  continuous in the present section.
We call semi-concave solution of (\ref{HJ}) a semi-concave function 
 $u:]0,T[\times \Rm^d\lto \Rm$ which solves (\ref{HJ}) at its points
of differentiability (these points form a set of full measure).
We recall that $D^e u(t,x)$ is the set of extreme points of the super-differential
$Du(t,x)$ and that it is contained in the set $D^*u(t,x)$ of reachable differentials,
see Section \ref{secnonsmooth}.

\begin{lem}
The semi-concave function $u(t,x)$ is a semi-concave solution if and only
if it satisfies the equality 
$\eta+H(t,x,p)=0$ for all $(t,x)$ and for all 
$(\eta,p)\in D^eu(t,x)$.
\end{lem}

\proof
If $u(t,x)$ is a solution, then the equality
$\partial_tu +H(t,x,\partial_xu)=0$ holds at each point of differentiability 
of $u$. Since $H$ is continuous we conclude that $\eta+H(t,x,p)=0$
for each reachable differential $(\eta,p)$ hence for each 
$(\eta,p)\in D^eu(t,x)$. 

The converse holds because, if $u$ is differentiable at $(t,x)$
then $(\partial_tu,\partial_xu)\in D^eu(t,x)$.
\qed

\begin{lem}
Let $u$ be a semi-concave solution of (\ref{HJ}).
Then $u$ is a viscosity solution if and only
each super-differential $(\eta,p)\in Du(t,x)$ 
satisfies the inequality $\eta+H(t,x,p)\leq 0$.
\end{lem}

\proof
A semi-concave solution is a viscosity super-solution.
Indeed, if $\phi(t,x)$ has a contact from below with $u$ at $(t_0,x_0)$,
then $u$ is differentiable at $(t_0,x_0)$, and
$$
\partial_t\phi(t_0,x_0)+H(t_0,x_0,\partial_x\phi(t_0,x_0))=
\partial_tu(t_0,x_0)+H(t_0,x_0,\partial_xu(t_0,x_0))=0.
$$
It is thus a solution if and only if it is a sub-solution.
\qed

The function $\check H^u_{t,x}$ is by definition the convex envelop of the function 
which is equal to $H_{t,x}$ on $D^e_xu(t,x)$ and equal to $+\infty$ at every other points.
It follows from \cite{Dac}, Theorem 2.35
that 
$$
\check H^u_{t,x}(p)=\inf_{p=\sum a_ip_i} \big(\sum a_iH(t,x,p_i)\big),
$$
where the infimum is taken on all the possibilities to write 
$p$ as the convex combination $\sum_{i=0}^{d}a_ip_i$ 
of $d+1$ points of $D^eu(t,x)$.
Similarly
$$
\hat H^u_{t,x}(p)=\sup_{p=\sum a_ip_i} \big(\sum a_iH(t,x,p_i)\big).
$$

Let us assume that  $u$ is a viscosity solution,
and fix a super-differential $(\eta,p)\in Du(t,x)$.
Since $Du(t,x)$ is the convex hull of $D^eu(t,x)$,
there exist $d+1$ positive numbers $a_i,1\leq 0\leq d$ and  points 
$(\eta_i,p_i)\in D^eu(t,x)$ such that 
$$
\sum_{i=0}^d a_i=1 \quad 
\text{and} \quad
\sum_{i=0}^d a_i(\eta_i,p_i)=(\eta,p).
$$  
Since $u$ is a semi-concave solution, each of the points 
$(\eta_i,p_i)$ solve the equation 
$
\eta_i+H(t,x,p_i)=0,
$
hence 
$$
\eta+\sum a_i H(t,x,p_i)=0.
$$
Since $u$ is a viscosity solution, we have the inequality
$$
\eta+H(t,x,p)\leq 0.
$$
We conclude that 
$$
H(t,x,p)\leq \sum a_i H(t,x,p_i)\leq \hat H^u_{t,x}(p).
$$
Since this holds for each point $(\eta, p)\in Du(t,x)$, we have
proved (\ref{concave}).

Let us now assume (\ref{convex}) and prove that 
$u$ is a viscosity solution.
We  fix a point $(t,x)\in ]0,T[\times \Rm^d$
and a point $(\eta,p)\in D^eu(t,x)$ and prove that 
$\eta+H(t,x,p)\leq 0$. 
Since $Du(t,x)$ is the convex hull of $D^eu(t,x)$, there exists
positive coefficients $a_i,0\leq i\leq d$ and points $(\eta_i,p_i)\in D^eu(t,x)$
such that $\sum_{i=0}^d a_i=1$ and $(\eta,p)=\sum_{i=0}^{d}a_i(\eta_i,p_i)$.
Using (\ref{convex}), we obtain
$$
\eta+H(t,x,p)
\leq \eta+\check H^u_{t,x}(t,x,p)\leq 
\sum_{i=0}^d a_i\big(\eta_i+ H(t,x,p_i)\big)=0
$$
since $\eta_i+H(t,x,p_i)=0$. This ends the proof of Theorem \ref{entropythm}.
\qed

\section{The Hopf formula}\label{secHF}
In this section, we
consider a continuous Hamiltonian which depends only on the variable $p$,
$H(t,x,p)=H(p)$, and  a (finite valued) concave initial condition 
$u_0$.
We write $u_0$ as an infimum  of affine
functions  with the formula
$
u_0(x)=\inf_p( p\cdot x -u_0^*(p)),
$
 where $u_0^*(p)$ is the (concave version of the) Legendre transform
$u_0^*(p):=\inf_x (p\cdot x-u_0(x))$. 
This dual  $u^*_0$  takes the value $-\infty$, outside of its  domain,
denoted by $P$.
Since  $u_0$ is finite valued  (hence locally bounded),
we have 
$$
\lim_{|p|\lto \infty } \frac{u^*_0(p)}{|p|}=-\infty.
$$
Since, in addition, the dual $u_0^*$ is upper semi-continuous, we have
$$
u_0(x)=\min_{p\in P}( p\cdot x -u_0^*(p)).
$$
The function $f(t,x)=px-u_0^*(p)-tH(p)$ is a
$C^2$ solution emanating from the affine initial 
condition $f_0(x)=px-u_0^*(p)$.
Let us define the function 
$$
u(t,x):=\inf_{p\in P} \big( px-u^*(p)-tH(p)\big)
$$
associated to this family.
Assume that  $T>0$ satisfies
$$
\lim_{|p|\lto \infty } \frac{u^*_0(p)+tH(p)}{|p|}=-\infty
$$
for each $t\in [0,T[$. This condition appears for example in \cite{LiR}.
In the case where $u_0$ is Lipschitz, $u_0^*$ is equal to $-\infty$
outside of a compact set and we can take $T= +\infty$.

The infimum defining $u$ is a minimum on $[0,T[\times \Rm^d$,
and we have (see Lemma \ref{surdiff})
$$
D^eu(t,x)\subset \{(-H(p),p) ,\quad p\cdot x-u_0^*(p)-tH(p)=u(t,p)\}
$$ 
for each $(t,x)\in [0,T[\times \Rm^d$, hence
$$
D^e_xu(t,x)\subset \{p\in \Rm^d ,\quad p\cdot x-u_0^*(p)-tH(p)=u(t,p)\}.
$$
We conclude that 
$$
H(p)=\frac{p\cdot x-u_0^*(p) -u(t,x)}{t}
$$
on $D^e_xu(t,x)$ and. Denoting by $c(p)$ the right hand side we have  $H(p)\leq c(p)$
on $\Rm^d$, hence on $D_xu(t,x)$.
Since $c(p)$ is a convex function of $p$,
this implies the entropy inequality (\ref{convex}).
We conclude from Theorem \ref{entropythm} that $u$ is a viscosity solution
on $[0,T[\times \Rm^d$.

In the case where $H$ also satisfies Hypothesis \ref{1}, and where $u_0$ is Lipschitz,
we conclude that $v=g=u$ on $[0,\infty)\times \Rm^d$.
The equality $g=u$ could also be deduced from Proposition \ref{hopf}.

\appendix
\section{A hessian computation}\label{hessian}
We consider the function
$f(x)=\varphi(|x|)$ on $\Rm^d$, where $\varphi:[0,\infty)\lto\Rm$ is a $C^2$ function.
Denoting by $e=x/|x|$ the radial direction, the Hessian of $f$ at a point $x\neq 0$ is
$$
d^2f_x[y,z]=
\frac{\varphi'(|x|)}{|x|}\langle y,z\rangle+
\frac{|x|\varphi''(|x|)-\varphi'(|x|)}{|x|}\langle e,y\rangle\langle e,z\rangle.
$$
In an orthonormal basis having the radial direction $e$ as first vector,
this bilinear form is expressed by a diagonal matrix 
with one diagonal  coefficient 
(corresponding to the radial direction)  equal to $\varphi''(|x|)$
and  $d-1$ diagonal coefficients  equal to $\varphi'(|x|)/|x|$.
As a consequence, its norm is 
$$
\|d^2f_x\|=\max \left\{|\varphi''(|x|)|, |\varphi'(|x|)|/|x|) \right\}.
$$
In the case considered in Section \ref{secnonsmooth},
we have $0\leq \varphi''(r)\leq B$ and $0\leq \varphi'(r)/r\leq B$
for each $r>0$, hence $\|d^2f_x\|\leq B$ at each point $x\neq 0$.
This inequality also obvioulsy holds at $x=0$ since 
$f(x)=B\|x\|^2/2$ near the origin.

\begin{footnotesize}

\end{footnotesize}
\end{document}